\documentclass{article}

\usepackage{amsmath, amsthm, amsfonts}

\newtheorem{thm}{Theorem}[section]

\newtheorem{lem}[thm]{Lemma}

\usepackage{mdframed}
\theoremstyle{definition}
\newmdtheoremenv{boxProb}{Problem}
\newmdtheoremenv{boxDef}{Definition}
\newmdtheoremenv{boxCor}{Corollary}
\newmdtheoremenv{boxThm}{Theorem}
\newmdtheoremenv{compjob}{Computational Job}
\newmdtheoremenv{reqi}{Requirement}

\usepackage{graphicx}
\usepackage{subfigure}          

\usepackage{placeins}

\usepackage{tabularx}

\usepackage{lmodern,textcomp}

\usepackage{algorithm}
\usepackage{algpseudocode}

\usepackage{enumitem}

\usepackage{xspace} 
\newcommand\largeparbreak{\par\bigskip}

\newcommand*\tageq{\refstepcounter{equation}\tag{\theequation}}

\newcommand{\inv}{^{-1}\xspace}
\newcommand{\bmu}{\boldsymbol{\mu}\xspace}

\newcommand{\hnu}{{\hat{\nu}}\xspace}
\newcommand{\blambda}{\boldsymbol{\lambda}\xspace}
\newcommand{\bxi}{\boldsymbol{\xi}\xspace}
\newcommand{\bs}{\textbf{s}\xspace}
\renewcommand{\t}{^\textsf{T}\xspace}

\newcommand{\hbz}{\hat{\textbf{z}}\xspace}

\newcommand{\away}[1]{}

\newcommand{\R}{\mathbb{R}\xspace}

\newcommand{\N}{\mathbb{N}\xspace}

\newcommand{\cA}{\mathcal{A}\xspace}
\newcommand{\cB}{\mathcal{B}\xspace}

\newcommand{\cL}{\mathcal{L}\xspace}

\newcommand{\cF}{\mathcal{F}\xspace}

\newcommand{\cO}{\mathcal{O}\xspace}

\newcommand{\bA}{\textbf{A}\xspace}

\newcommand{\bG}{\textbf{G}\xspace}

\newcommand{\tbH}{\widetilde{\textbf{H}}\xspace}
\newcommand{\tbM}{\widetilde{\textbf{M}}\xspace}

\newcommand{\bD}{\textbf{D}\xspace}
\newcommand{\bS}{\textbf{S}\xspace}

\newcommand{\bQ}{\textbf{Q}\xspace}
\newcommand{\bK}{\textbf{K}\xspace}

\newcommand{\bx}{\textbf{x}\xspace}

\newcommand{\bb}{\textbf{b}\xspace}

\newcommand{\bc}{\textbf{c}\xspace}

\newcommand{\bd}{\textbf{d}\xspace}

\newcommand{\hbx}{\hat{\textbf{x}}\xspace}

\newcommand{\bz}{\textbf{z}\xspace}
\newcommand{\be}{\textbf{1}\xspace}
\newcommand{\bei}[1]{{\textbf{e}}\xspace}

\newcommand{\bM}{\textbf{M}\xspace}

\newcommand{\bH}{\textbf{H}\xspace}
\newcommand{\bI}{\textbf{I}\xspace}

\newcommand{\bO}{\textbf{0}\xspace}

\newcommand{\tbx}{\tilde{\textbf{x}}\xspace}

\newcommand{\bv}{\textbf{v}\xspace}

\newcommand{\opdiag}{\textsl{diag}\xspace}

\newcommand{\cond}{\textsl{cond}\xspace}

\newcommand{\tol}{{\textsf{tol}}\xspace}

\newcommand{\epsMach}{\varepsilon_{\textsf{mach}}\xspace}

\newcommand{\cbx}{\check{\bx}\xspace}
\newcommand{\calpha}{\check{\alpha}\xspace}

\newcommand{\Ipopt}{\textsc{Ipopt}\xspace}

\title{PPD-IPM: Outer primal, inner primal-dual interior-point method for nonlinear programming}

\author{Martin Neuenhofen}

\date{\today}

\begin{document}

\maketitle

\begin{abstract}
	In this paper we present a novel numerical method for computing local minimizers of twice smooth differentiable non-linear programming (NLP) problems.
	
	So far all algorithms for NLP are based on either of the following three principles: successive quadratic programming (SQP), active sets (AS), or interior-point methods (IPM). Each of them has drawbacks. These are in order: iteration complexity, feasibility management in the sub-program, and utility of initial guesses. Our novel approach attempts to overcome these drawbacks.
	
	We provide: a mathematical description of the method; proof of global convergence; proof of second order local convergence; an implementation in \textsc{Matlab}; experimental results for large sparse NLPs from direct transcription of path-constrained optimal control problems.
\end{abstract}


\section{Introduction}
We consider the numerical solution of the following non-linear programming (NLP) problem for local minimizers:
\begin{subequations}
\begin{align}
	&\min_{\bx \in \R^n} 	& 	f(\bx)& 		\\
	&\text{s.t.} 			&	c(\bx)&=\bO\,,	\\
	& 						& -\be\leq\phantom{g(}\bx\phantom{)}&\leq \be
\end{align}\label{eqn:NLP}
\end{subequations}
where $\be \in \R^n$ is the vector of ones and
\begin{align*}
	f\,&:\, \R^n \rightarrow \R^1\,, & c\,&:\, \R^n \rightarrow \R^m
\end{align*}
are bounded twice Lipschitz-continuously differentiable functions. We write $\bx^\star$ for an arbitrary fixed local minimizer. The Lagrange-function is defined as
\begin{align*}
	\cL\,:\,\R^n \times \R^m \rightarrow \R^1, (\bx,\blambda) \mapsto f(\bx) - \blambda\t\cdot c(\bx)\,.
\end{align*}
Every non-linear programming problem can be substituted into form \eqref{eqn:NLP} by using bounds on $\|\bx^\star\|_\infty$ and slacks for inequality constraints. The dimensions are $m,n \in \N$, where $m$ can be smaller, equal or larger than $n$.

\paragraph{Penalty barrier program} In this paper we treat \eqref{eqn:NLP} by solving a related minimization problem. The Karush-Kuhn-Tucker (KKT) conditions for \eqref{eqn:NLP} are
\begin{subequations}
\begin{align}
\nabla f(\bx) + \rho \cdot \bS \cdot \bx - \nabla c(\bx)\cdot \blambda - \bmu_L + \bmu_R = &\bO \label{eqn:KKT:NLP:1}\\
c(\bx) + \omega \cdot \blambda = &\bO \label{eqn:KKT:NLP:2}\\
\opdiag(\bmu_L) \cdot (\be+\bx) - \tau \cdot\be = &\bO \label{eqn:KKT:NLP:3}\\
\opdiag(\bmu_R) \cdot (\be-\bx) - \tau \cdot\be = &\bO \label{eqn:KKT:NLP:4}\\
\bmu_L \geq& \bO \label{eqn:KKT:NLP:5}\\
\bmu_R \geq& \bO \label{eqn:KKT:NLP:6}\\
\be+\bx \geq& \bO \label{eqn:KKT:NLP:7}\\
\be-\bx \geq& \bO\,,\label{eqn:KKT:NLP:8}
\end{align}\label{eqn:KKT:NLP}
\end{subequations}
where $\rho=\omega=\tau=0$, $\blambda \in \R^m$, $\bmu_L,\bmu_R \in \R^n_+$, and $\bS \in \R^{n \times n}$ symmetric positive definite. These equations can be numerically unsuitable. E.g., when $\nabla f(\bx)$ and $c(\bx)$ have constant values around some $\bx$ then the system can be locally non-unique for $\bx$. And when columns of $\nabla c(\bx)$ are linearly dependent then there are multiple solutions for the Lagrange multiplier $\blambda$ at fixed $\bx$. Also $\bmu_L,\bmu_R$ can be non-unique when inequality constraints are active whose gradients are linearly dependent to the columns of $\nabla c(\bx)$.

For sufficiently small regularization parameters $\rho,\omega,\tau>0$ the system's solution is locally unique, which is desirable when using equations \eqref{eqn:KKT:NLP:1}--\eqref{eqn:KKT:NLP:4} within Newton's iteration to compute local solutions $\bx,\blambda,\bmu_L,\bmu_R$. This is because the uniqueness gives regularity of the Jacobian that appears in the linear equation system of Newton's method. Thus, there is second order local convergence of the iterates. Substituting \eqref{eqn:KKT:NLP:2}--\eqref{eqn:KKT:NLP:4} into \eqref{eqn:KKT:NLP:1}, we find that solutions of \eqref{eqn:KKT:NLP} are critical points for the following problem.
\begin{align}
	\min_{\bx \in \overline{\Omega}} \quad \phi(\bx) \label{eqn:MinPhi}
\end{align}
where
\begin{align*}
	\phi(\bx)&:=f(\bx) + \frac{\rho}{2} \cdot \|\bx\|_\bS^2 + \frac{1}{2 \cdot \omega} \cdot \|c(\bx)\|_2^2 - \tau \cdot \be\t\cdot\Big(\, \log(\be+\bx) + \log(\be-\bx) \,\Big) \\[3pt]
	\Omega &:= \big\lbrace \bxi \in \R^n \ : \ \|\bxi\|_\infty < 1 \big\rbrace\,,
\end{align*}
$\|\bx\|_\bS := \sqrt{\bx\t\cdot\bS\cdot\bx}$ is the induced norm and $\log(\cdot)$ is the natural logarithm of each component of the argument. The equations \eqref{eqn:KKT:NLP:5}--\eqref{eqn:KKT:NLP:8} are strictly forced because $\phi$ goes to infinity as $\bx$ approaches $\partial\Omega$. Numerically suitable values for the regularization are $\rho,\omega,\tau$ between $10^{-5}$ and $10^{-8}$.

The program \eqref{eqn:MinPhi} we call \textit{penalty-barrier program}. It is badly scaled for small values of $\rho,\omega,\tau>0$. This is why iterative schemes based on \mbox{(Quasi-)}Newton-type descent directions yield poor progress for it and would result in an impractically large amount of iterations \cite[p. 569ff, p. 621]{Boyd}, \cite{SUMT}.

\paragraph{Outer primal inner primal-dual method}
In this paper we present a novel approach to solving \eqref{eqn:MinPhi}. We still perform a direct minimization of \eqref{eqn:MinPhi} because we believe that this is the robustest approach. Since search directions from Newton steps would yield bad progress, we instead use search directions that are obtained from the solution of subproblems of the following form:
\begin{align}
\min_{\bx \in \overline{\Omega}} \quad q(\bx) \label{eqn:MinQ}
\end{align}
where
\begin{align*}
q(\bx):=&\frac{1}{2}\cdot\bx\t\cdot\bQ\cdot\bx+\bc\t\cdot\bx + \frac{1}{2 \cdot \omega} \cdot \|\bA\cdot\bx-\bb\|_2^2 \\
&\quad - \tau \cdot \be\t\cdot\Big(\, \log(\be+\bx) + \log(\be-\bx) \,\Big)\tageq\label{eqn:def:q}
\end{align*}
and $\bQ \in \R^{n \times n}$ is symmetric positive definite, $\bc \in \R^n$, $\bA \in \R^{m \times n}$, $\bb \in \R^m$. In our algorithm $\bQ$ will be (approximately)
\begin{align*}
	\nabla^2_{\bx\bx} \cL(\bx,\blambda) + \rho \cdot \bS
\end{align*}
so that at a given iterate $\bx,\blambda$ it holds $\nabla q = \nabla \phi$ and (approximately) $\nabla^2 q = \nabla^2 \phi$.

Problems \eqref{eqn:MinQ} can be solved efficiently using a particular primal-dual path-following method described in \cite{StableIPM}. The avoidance of a quadratic approximation for the logarithmic terms yields a better fit of the search direction to minimize \eqref{eqn:MinPhi}. Since we minimize \eqref{eqn:MinPhi} directly with the search directions obtained from \eqref{eqn:MinQ}, there is no need to spend extra attention on the convergence of feasibility: The equality constraints $c(\bx)=\bO$ are treated with a quadratic penalty that is well-represented in $q$. We will employ a watch-dog technique \cite{watchdog} to achieve large steps along the directions obtained from \eqref{eqn:MinQ} even though the penalty parameter $\omega>0$ is very small. The inequalities $-\be\leq\bx\leq\be$ are forced through barriers in \eqref{eqn:MinPhi}. These are considered in an unmodified way in \eqref{eqn:MinQ}, always keeping $\bx \in \Omega$. Altogether this results in a simple and robust algorithm that is easy to implement and analyse.

\subsection{Literature review}
Algorithms for NLP can be divided into three distinct classes, confer to \cite{NumOpt}: active set methods (ASM), successive quadratic programming (SQP), and interior-point methods (IPM).
\largeparbreak

ASM are based on iteratively improving a guess of the active inequality constraints in \eqref{eqn:NLP}. The guess is stored as a set $\cA$ of indices, called \textit{active set}. Using a guess for $\cA$, an equality-constrained non-linear programming problem is formed and solved for a local minimizer $\bx^\star_\cA$. At $\bx^\star_\cA$ the Lagrange multipliers provide information on the optimality of $\cA$. If $\cA$ is non-optimal, then a new estimate $\cA$ for the active set is formed and again $\bx^\star_\cA$ is computed. This procedure is repeated until $\cA$ is correct, which implies $\bx^\star_\cA \equiv \bx^\star$ is a local minimizer of \eqref{eqn:NLP}. An introduction to active set methods can be found in \cite{Flet87}.

ASM are numerically robust because no penalty or barrier must be introduced to treat the constraints. As a further advantage, ASM provide additional information on the set $\cA$ of active constraints at the local minimizer. The problem however with ASM is that there is no polynomially efficient method for determining the optimal active set $\cA$. Problems are known for which ASM would try all possible active sets \cite{klee1970good} until in the very last attempt they find the correct one. This results in a worst-case time complexity that grows exponentially with $n$ \cite{klee1970good}.

\largeparbreak
SQP methods improve the current iterate by moving in a direction obtained from solving a convex quadratic sub-program. The step-size along this direction is determined by minimizing a merit-function or using a filter. For a general overview on SQP methods consult \cite{Gill:2005:SSA:1055334.1055404}.

Special care must be taken to modify the sub-program accordingly such that it always admits a feasible solution. Typically this is achieved through $\ell_1$-penalties. This is sometimes referred to as \textit{elastic mode} \cite{SQPinconsistent}. The $\ell_1$-penalties in the subproblem must be sufficiently large to ensure progress towards feasibility. On the other hand, too large values for the penalties lead to a bad scaling of the quadratic subproblem, confer to \cite{Han1977}.

Special care must be further taken to make sure that ---despite the modification with the elastic mode--- the search direction obtained from the subproblem is still a descent direction for the line-search. As one possible way to achieve this, the penalty parameters in the $\ell_1$-merit-function must be chosen with respect to those in the subproblem, confer to \cite{Schittkowski1982}. If the penalty terms in the merit-function are too large then it is likely that the line-search admits small steps only, confer to \cite{watchdog}. This can be resolved, e.g., by using second-order corrections \cite{TrustRegionMethods}, which however may require the computationally prohibitive task of solving a convex quadratic program at several trial points.

The sub-program that must be solved in each iteration is a convex quadratic program (CQP). CQP can be solved using either active set methods or interior-point methods. Active set methods can have exponential time complexity in the worst case but can be fast in practice. In contrast to that, there are interior-point methods for CQP that are proven to be polynomially efficient in theory \cite{Kar84}. In practice they converge very fast. The field is strongly influenced by Mehrotra's predictor-corrector method \cite{Mehrotra}, which is a primal-dual interior-point method that can be used for solving CQP in a very efficient way.
\largeparbreak

IPM solve \eqref{eqn:NLP} by considering a barrier function as in \eqref{eqn:MinPhi}. The inequality constraints are removed and instead the cost-function is augmented with so-called barrier terms. These are terms that go to infinity when $\bx$ approaches the border of $\Omega$. The barrier-augmented cost-function we call $f_\tau$. For example, $f_\tau$ could be
\begin{align*}
	f_\tau(\bx) = f(\bx)- \tau \cdot \be\t\cdot\Big(\, \log(\be+\bx) + \log(\be-\bx) \,\Big)\,.
\end{align*}
For small values $\tau>0$, e.g. $\tau=10^{-10}$, the barrier-term mildly influences the level-sets of $f_\tau$ in the interior of $\Omega$. All minimizers of $f_\tau$ are interior and thus satisfy the inequality constraints in a strict way.

To make sure that the unconstrained minimizers of $f_\tau$ are accurate approximations to the constrained minimizers of $f$ it is necessary to choose $\tau>0$ very small. However, for small $\tau$ the barrier term leads to a bad scaling of the barrier-augmented cost function. This results in bad progress when using descent directions obtained from \mbox{Quasi-}Newton-type methods, which however are used in almost every IPM, compare e.g. to \cite{IPOPT,LOQO,HOPDM}. This is why practical algorithms decrease the size of $\tau$ iteratively within the iterative computation of $\bx$. Thus, initially $\tau$ is large and yields good progress for the iterates of $\bx$. As $\bx$ approaches the minimizer, $\tau$ is slowly reduced and $\bx$ needs only be mildly refined. For an introduction to interior-point methods we refer to \cite{Wright,IPM25ylater}.

For many classes within the domain of convex programming there is strong evidence on the computational efficiency of IPM. Prominent examples are primal methods for self-concordant functions \cite{nesterov1994interior} and primal-dual methods for linear programming \cite{Wright}. However, for general NLP there is no result available on the complexity of the iteration count of IPM; compare to \cite{Forsgren02interiormethods}.

A serious disadvantage of IPM it their difficulty to make good use of initial guesses \cite{YildirimW02}. This phenomenon can be explained by the fact that for the initially large values of $\tau$ the function $f_\tau$ has little in common with $f$. Thus, a potentially good initial guess $\bx_0$ of the local minimizer is driven away in early iterations of IPM towards a minimizer of $f_\tau$ for this initially large value of $\tau$. Eventually, $\tau$ decreases and the iterates $\bx$ move back the local minimizer (to which the initial guess may have been close, or to another one).

\largeparbreak
In contrast to interior-point methods we use a fixed value of $\tau$ within the minimization of \eqref{eqn:MinPhi}. Thus, our method does not move away from good initial guesses if they are close to local minimizers of $\phi$. A strategy with decreasing values for $\tau$ is not required in our method because even for small values like $\tau=10^{-8}$ the search directions obtained from solving \eqref{eqn:MinQ} allow fast progress within the line-search on $\phi$. This holds because the value of $\tau$ does not influence the accuracy in which $q$ approximates $\phi$.

\subsection{Structure}
In Section~2 we present the numerical method. In Section~3 we provide proofs for the convergence: We prove that the local minimizers of \eqref{eqn:MinPhi} converge to the constrained local minimizers of \eqref{eqn:NLP}. We prove global convergence of our numerical method and we prove second order local convergence. In Section~4 we discuss details of our implementation and practical enhancements. Section~5 presents numerical experiments against \Ipopt \cite{IPOPT} for large sparse non-linear programs that arise from the direct discretization of optimal control problems. Eventually we draw a summarizing conclusion.

\section{Primal-primal-dual interior-point method}
Our method is an iterative method that computes a sequence $\lbrace \bx_k \rbrace \subset \Omega$ of interior-points that converge to stationary points of \eqref{eqn:MinPhi}. $\rho,\omega,\tau>0$ and $\bS\in\R^{n \times n}$ symmetric positive definite are considered to be provided by the user. The values $\rho=\omega=\tau=10^{-7}$ and $\bS = \bI_{n \times n}$ are often suitable.

The method goes as follows. Given $\bx_k$, either from a former iteration or an initial guess when $k=0$, we compute
\begin{align*}
	\blambda_k := \frac{-1}{\omega}\cdot c(\bx_k)\,.
\end{align*}
We notice
\begin{align*}
\nabla \phi(\bx_k)&= \rho \cdot \bS \cdot \bx_k + \nabla_\bx \cL(\bx_k,\blambda_k) - \frac{\tau}{\be+\bx_k} + \frac{\tau}{\be-\bx_k}\\
\nabla^2 \phi(\bx_k)&= \rho \cdot \bS + \bM_k + {\tau}\cdot\Big( \opdiag(\be+\bx_k)^{-2} + \opdiag(\be-\bx_k)^{-2}\Big)
\end{align*}
where
\begin{subequations}
\begin{align}
	\bM_k &:= \bH_k + \frac{1}{\omega}\cdot\nabla c(\bx_k) \cdot \nabla c(\bx_k)\t\\
	\bH_k &:= \nabla^2_{\bx\bx} \cL(\bx_k,\blambda_k)\,.
\end{align}\label{eqn:def:MH}
\end{subequations}
If 
\begin{align}
	\nabla_{\bx\bx}^2\cL(\bx_k,\blambda_k)
\end{align}
is symmetric positive semi-definite then $\bH_k$ and thus $\bM_k$ are positive semi-definite. If
\begin{align}
	\big(\nabla c(\bx_k)^\perp\big)\t \cdot \nabla_{\bx\bx}^2\cL(\bx_k,\blambda_k) \cdot \big(\nabla c(\bx_k)^\perp\big)\label{eqn:projectedhessian}
\end{align}
is positive semi-definite ---\,which is a necessary condition at least for interior minimizers $\bx^\star$ of \eqref{eqn:NLP}\,--- then there exist suitably small values $\omega>0$ such that $\bM_k$ is positive semi-definite. Whenever $\bM_k$ is positive semi-definite it follows in turn that $\nabla^2 \phi(\bx_k)$ is positive definite.

We form a \textit{penalty-barrier convex quadratic program} \eqref{eqn:MinQ} with the following values
\begin{subequations}
\begin{align}
	\bQ_k &:= \tbH_k+\rho\cdot\bS\\
	\bA_k &:= \nabla c(\bx_k)\t\\
	\bc_k &:= \nabla_\bx f(\bx_k) \ - \tbH_k \cdot \bx_k\\
	\bb_k &:= - c(\bx_k) \ + \bA_k \cdot \bx_k
\end{align}\label{eqn:def:cqp}
\end{subequations}
where $\tbH_k$ is an approximation to $\bH_k$ such that
$$\tbM_k := \tbH_k + \frac{1}{\omega}\cdot\nabla c(\bx_k) \cdot \nabla c(\bx_k)\t$$ 
is symmetric positive semi-definite. By construction, the resulting function $q$ in \eqref{eqn:def:q} satisfies
\begin{align*}
	\nabla q(\bx_k) &= \nabla \phi(\bx_k)
\end{align*}
and $q$ is strictly convex due to positive semi-definiteness of $\tbM_k$. If in addition one of the above-mentioned conditions holds then the choice $\tbH_k = \bH_k$ is suitable such that $\tbM_k=\bM_k$ is symmetric positive semi-definite. It then follows
\begin{align*}
\nabla^2 q(\bx_k) &= \nabla^2 \phi(\bx_k)\,.
\end{align*}

Now that $q$ is fully defined, we solve the penalty-barrier convex quadratic program \eqref{eqn:MinQ} of it. In \cite{StableIPM} we describe a short-step primal-dual path-following method that can solve \eqref{eqn:MinQ} in weakly polynomial time complexity. The method described in the reference is further numerically stable if $\epsMach$ is chosen sufficiently small with respect to a weak constant that depends on the logarithms of the norms of $\bQ_k,\bA_k,\bc_k,\bb_k$ and the logarithm of $\omega$. In Section~4 of this paper we provide a long-step variant of the referred method that is fast and reliable in practice, is suitable also for large sparse problems, and can solve \eqref{eqn:MinQ} to high numerical accuracy.

Once that the solution of \eqref{eqn:MinQ} is obtained, we write it into a vector $\hbx_k$. We define the step-direction $\bv_k := \hbx_k - \bx_k$. Finally, we compute a new iterate
\begin{align*}
	\bx_{k+1} := \bx_k + \alpha_k \cdot \bv_k
\end{align*}
where $\alpha_k \in\R^+$ is chosen to minimize $\phi$ along the line $\bx(\alpha) := \bx_k + \alpha \cdot \bv_k \in \Omega$. $\alpha$ is chosen, e.g., by using a back-tracking line-search with Armijo-rule, confer to \cite{NumOpt}. Algorithm~\ref{algo:Solver} encapsulates the algorithmic steps.

\begin{algorithm}
	\caption{PPD-IPM, pure version}
	\label{algo:Solver}
	\begin{algorithmic}[1]
		\Procedure{PPDIPM}{\,$\bx_0,\rho,\omega,\tau,\bS,\tol$\,}
		\State $k:=0$
		\While{$ \|\nabla \phi(\bx_k)\|_2> \tol $}
			\State $\blambda_k := \frac{-1}{\omega}\cdot c(\bx_k)$
			\State Choose $\tbH_k \approx \bH_k$ such that $\tbM_k$ is positive semi-definite.
			\State Compute $\bQ_k,\bc_k,\bA_k,\bb_k$ from \eqref{eqn:def:cqp}, defining $q$.
			\State Compute $\hbx_k$, the minimizer of \eqref{eqn:MinQ}.
			\State $\bv_k := \hbx_k - \bx_k$
			\State $\alpha_k := \operatornamewithlimits{argmin}_{\alpha \in \R^+}\big\lbrace \phi(\bx_k + \alpha \cdot \bv_k) \big\rbrace$ \label{algo:Solver:lineAlpha}
			\State $\bx_{k+1} := \bx_k + \alpha_k \cdot \bv_k$
			\State $k := k+1$\label{algo:Solver:line:End}
		\EndWhile
		\State \Return $\bx_k$
		\EndProcedure
	\end{algorithmic}
\end{algorithm}

\section{Proof of convergence}
We show that the iterates $\bx_k$ of Algorithm~\ref{algo:Solver} converge to stationary points of \eqref{eqn:MinPhi}. In the first subsection we show that there is convergence to a stationary point of $\phi$ from every initial guess $\bx_0$. In the second subsection we show that under suitable conditions there is second-order convergence of $\|\nabla\phi(\bx_k)\|_2$ to zero for $k \in \N$ greater than some finite number.

\subsection{Global convergence}

We start with some technical results.

\begin{lem}[Boundedness]\label{lem:BoundednessOmega0}
	Let $\bx \in \Omega$. We define
	\begin{align*}
		\overline{\Omega}_0(\bx) := \lbrace \bxi \in \Omega \ : \ \phi(\bxi) \leq \phi(\bx)\, \rbrace\,.
	\end{align*}
	The space $\overline{\Omega}_0(\bx)$ is always bounded and closed. Further, for each $\overline{\Omega}_0(\bx)$ there is a constant $C_\phi\in \R$ such that
	\begin{align*}
		|\phi(\tbx)|,\ \|\nabla\phi(\tbx)\|_2,\ \|\nabla^2 \phi(\tbx)\|_2 \leq C_\phi \quad \forall \tbx \in \overline{\Omega}_0(\bx)\,.
	\end{align*}
\end{lem}
\begin{itshape}
	\noindent 
	\underline{Proof:} left to the reader. q.e.d.
\end{itshape}

\begin{lem}[Sufficient descent]\label{lem:SufficientDescent}
	Let $\bx \in \Omega$, $\bv \in \R^n \setminus \lbrace \bO \rbrace$, $0 < \vartheta < \pi/2$. If the angular condition
	\begin{align*}
		\angle(\,\bv\,,\,-\nabla \phi(\bx)\,)\leq \frac{\pi}{2} - \vartheta
	\end{align*}
	is satisfied then $\exists \ \theta > 0$, only depending on $\vartheta$ and $\|\nabla\phi(\bx)\|_2$, such that the following holds:
	\begin{align*}
		\min_{\alpha \in \R^+} \Big\lbrace\, \phi(\bx + \alpha \cdot \bv) \,\Big\rbrace \leq \phi(\bx)-\theta
	\end{align*}
\end{lem}
\begin{itshape}
\noindent 
\underline{Proof:} Confer to \cite[Sections~3.1--3.2]{NumOpt}.
\end{itshape}

\begin{figure}
	\centering
	\includegraphics[width=0.999\linewidth]{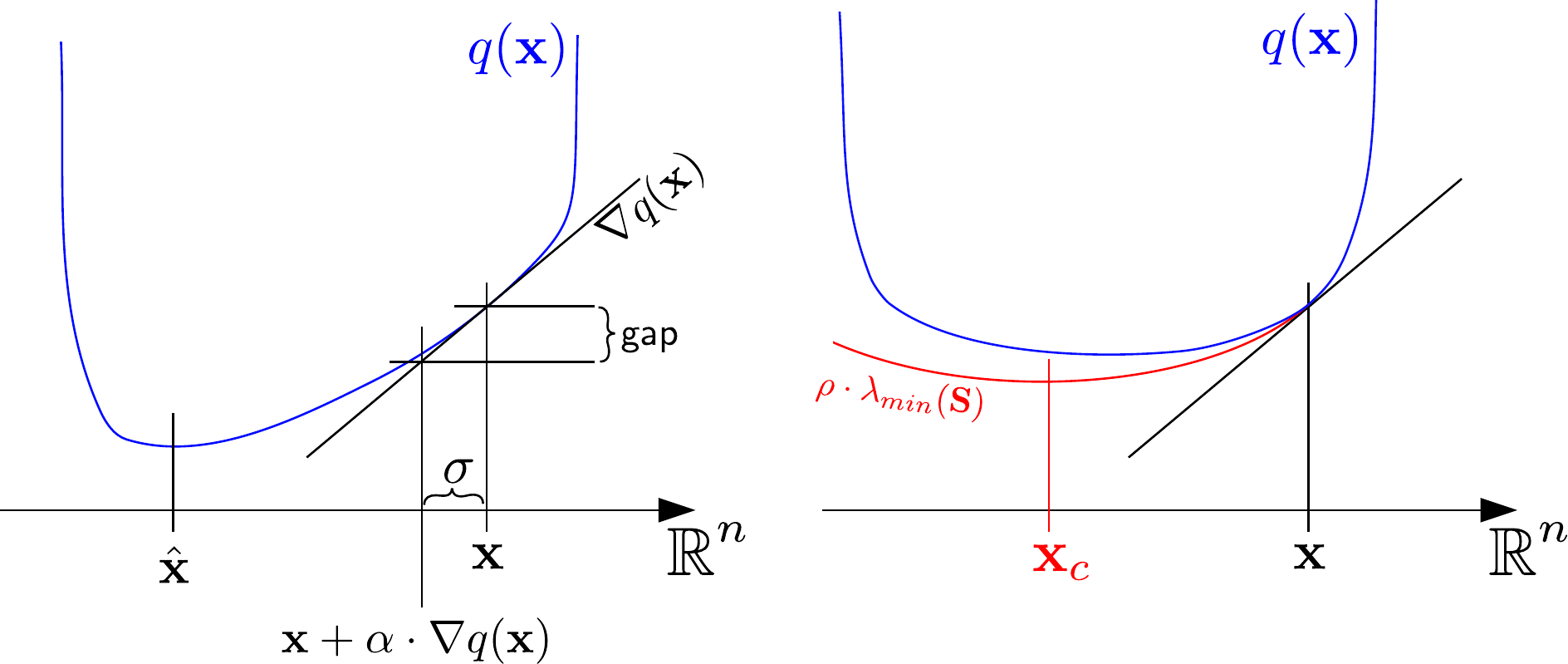}
	\caption{Left: Plot of $q$ through line $\hbx \rightarrow \bx$. $q$ is convex and slope from $\bx$ to $\hbx$ is bounded below by $\|\nabla q(\bx)\|_2$, thus $\sigma$ bounds distance $\|\hbx-\bx\|_2$ from below. Right: Red curve is quadratic function $\psi$ with isotropic second derivative $\rho \cdot \lambda_\text{min}(\bS)$. Thus $q$ is bounded below by $\psi$. We can bound the distance of $\hbx$ to $\bx_c$.}
	\label{fig:proofdescentbsmall}
\end{figure}

\begin{figure}
	\centering
	\includegraphics[width=0.99999\linewidth]{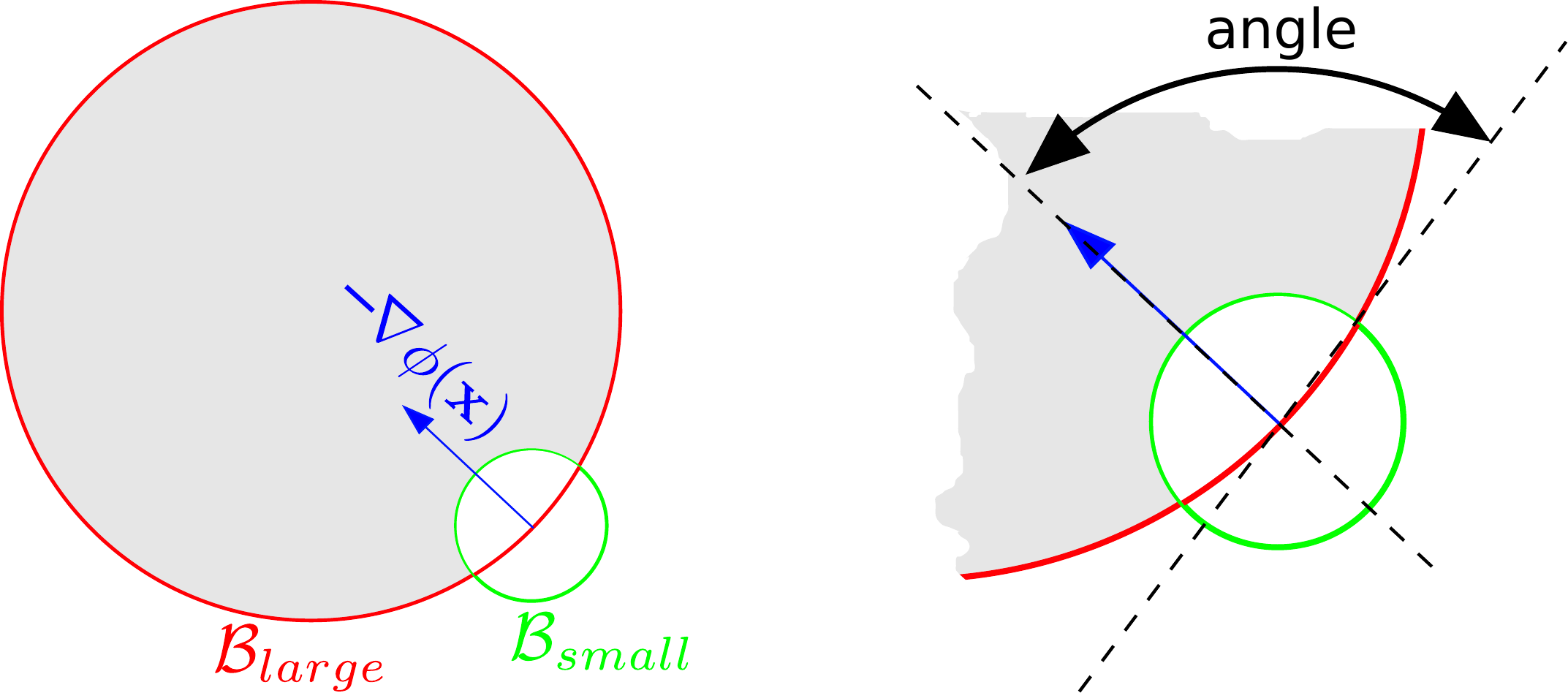}
	\caption{Plot of $\cB_{\text{large}}$ and $\cB_{\text{small}}$. They only depend on $\nabla\phi(\bx)$ and $C_{\phi}$, where the latter only depends on $\overline{\Omega}(\bx_0)$ for some former initial guess from which $\bx$ may have propagated. $\hbx$ lives in the gray region, implying that $\hbx-\bx$,$-\nabla\phi(\bx)$ have an angle of strictly less than $90$ degrees.}
	\label{fig:blargebsmall}
\end{figure}

\begin{lem}[Sufficient descent direction]\label{lem:SufficientDescentDirection}
	Let $\rho>0$, $\lambda_\text{min}(\bS)>0$, $\bx \in \Omega$. Then the following holds:
	\begin{align*}
		\forall \eta > 0 \quad \exists \vartheta > 0 \quad : \quad \|\nabla \phi(\bx)\|_2 \geq \eta \ \Rightarrow \ \angle(\,\hbx-\bx\,,\,-\nabla\phi(\bx)\,) \leq \frac{\pi}{2}-\vartheta
	\end{align*}
\end{lem}
\begin{itshape}
	\noindent 
	\underline{Proof:} Consider $\hbx$, computed as local minimizer of the local approximation function $q(\cdot)$ of $\phi(\cdot)$ around $\bx$, defined in \eqref{eqn:MinQ}. Our proof works by showing that $\hbx \in \cB_{\text{large}} \setminus \cB_{\text{small}}$ holds, where $\cB_{\text{large}}\,,\ \cB_{\text{small}}$ are two spheres. The geometric relation of the two spheres then enforces the claimed angular condition.
	
	We start with $\cB_{\text{small}}$. From the steepest descent direction at $\bx$ we find
	\begin{align*}
		\min_{\tbx \in \Omega}\big\lbrace\,q(\tbx)\,\big\rbrace \leq q\big(\bx + \alpha \cdot \overbrace{\nabla q(\bx)}^{\equiv \nabla \phi(\bx)}\big) \leq q(\bx) - \alpha \cdot \|\nabla \phi(\bx)\|_2 + \frac{\alpha^2}{2} \cdot \|\nabla^2\phi(\bx)\|_2^2 \cdot C_{\phi}\,.
	\end{align*}
	Inserting $\alpha = \frac{1}{C_{H\phi}}$, we get
	\begin{align*}
		\min_{\tbx \in \Omega}\big\lbrace\,q(\tbx)\,\big\rbrace \leq q(\bx) - \underbrace{\frac{\|\nabla\phi(\bx)\|_2^2}{2 \cdot C_{\phi}}}_{=:\textsf{gap}}
	\end{align*}
	Since $\hbx = \operatornamewithlimits{argmin}_{\tbx \in \Omega}\lbrace q(\tbx)\rbrace$, and since the negative slope of $q$ below $q(\bx)$ is bounded by $\|\nabla\phi(\bx)\|_2$, cf. Figure~\ref{fig:proofdescentbsmall} left, we find
	\begin{align*}
		\|\hbx-\bx\|_2 \geq \sigma = \frac{\textsf{gap}}{\|\nabla\phi(\bx)\|_2} = \frac{\|\nabla\phi(\bx)\|_2}{C_{\phi}}\,.
	\end{align*}
	We define $\cB_{\text{small}}:= \lbrace \,\bxi \in \R^n\ : \ \|\bxi-\bx\|_2<\sigma\,\rbrace$.
	
	Consider 
	\begin{align*}
		\psi(\tbx) = \phi(\bx) + \nabla\phi(\bx)\t \cdot \tbx + \frac{\rho\cdot\lambda_\text{min}(\bS)}{2}\cdot \|\tbx\|_2^2\,
	\end{align*}
	cf. Figure~\ref{fig:proofdescentbsmall} right. We define the minimizer
	\begin{align*}
		\bx_c := \bx - \frac{1}{\rho \cdot \lambda_\text{min}(\bS)} \cdot \nabla\phi(\bx)
	\end{align*}
	of $\psi$ and $\cB_{\text{large}} := \lbrace \bxi \in \R^n \ : \ \|\bxi - \bx_c\|_2 \leq \|\bx-\bx_c\|_2 \rbrace$. Since $\psi(\cdot)$ is a lower bound on $q(\cdot)$ it must hold
	$\hbx \in \cB_{\text{large}} \setminus \cB_{\text{small}}$. Now consider Figure~\ref{fig:blargebsmall}, from which we find that the claimed angular condition must hold. q.e.d.
\end{itshape}

\begin{lem}[Limit point]\label{lem:LimitPoint}
	Choose the initial guess $\bx_0 \in \Omega$ and consider the iterates $\lbrace \bx_k \rbrace_{k \in \N_0} \subset \overline{\Omega}(\bx_0)$ of Algorithm~\ref{algo:Solver}. Define $\bx_\infty := \lim\limits_{k \rightarrow \infty} \bx_k$. Then $\forall \varepsilon>0$ $\exists N \in \N$ such that the following holds $\forall k \geq N$:
	\begin{align*}
		 |\phantom\nabla\phi(\bx_\infty) - \phantom\nabla\phi(\bx_k) |\phantom{\|_2} \leq & \varepsilon\\
		\|\nabla\phi(\bx_\infty) - \nabla\phi(\bx_k) \phantom|\|_2 \leq & \varepsilon
	\end{align*}
\end{lem}
\begin{itshape}
\noindent\underline{Proof:} follows by L-continuity and Lemma~\ref{lem:BoundednessOmega0}. q.e.d.
\end{itshape}

Now we have everything in hand for the final result. The following theorem proves the global convergence of Algorithm~\ref{algo:Solver} to a stationary point of \eqref{eqn:MinPhi}.
\begin{thm}[Stationary limit]\label{lem:StationaryLimit}
	Consider the properties from Lemma~\ref{lem:LimitPoint}. Then:
	\begin{align*}
		\|\nabla\phi(\bx_\infty)\|_2 = 0
	\end{align*}
\end{thm}
\begin{itshape}
\noindent\underline{Proof:} (by contradiction). Let $\bd := \nabla\phi(\bx_\infty)$, $\eta := 0.5 \cdot \|\bd\|_2$.
\begin{align}
	\text{We assume }\eta > 0\,.\label{eqn:lem:StationaryLimit:WrongAssumption}
\end{align}
We choose a strictly positive value $\varepsilon < \eta$ for Lemma~\ref{lem:LimitPoint} and get $N \in \N$. Consider an arbitrary integer $k \geq N$. Notice that from Lemma~\ref{lem:LimitPoint} follows
\begin{align*}
	\phi(\bx_\infty)-\varepsilon \leq \phi(\bx_{k+1})\,.\tageq\label{eqn:StationaryLimitContra1}
\end{align*}
Since $\|\nabla\phi(\bx_k)-\bd\|_2\leq \eta$ holds according to Lemma~\ref{lem:LimitPoint}, it follows $\|\nabla\phi(\bx_k)\|_2\geq \eta$. We apply Lemma~\ref{lem:SufficientDescentDirection} to obtain $\vartheta>0$. Notice that $\vartheta,\eta$ are independent of $k$.

Define $\bv_k := \\hbx_k - \bx_k$. Due to Lemma~\ref{lem:SufficientDescentDirection}, $\bx_k$, $\bv_k$ and $\vartheta$ together satisfy the angular condition of Lemma~\ref{lem:SufficientDescent}, which then says that $\bx_{k+1}=\bx_k + \alpha_k \cdot \bv_k$ satisfies
\begin{align*}
	\phi(\bx_{k+1}) \leq \phi(\bx_k) - \theta\tageq\label{eqn:StationaryLimitContra2}
\end{align*}
where $\theta$ depends only on $\vartheta,\eta$, which in turn do not depend on $k,\varepsilon$. Thus, we can choose $\varepsilon>0$ sufficiently small such that \eqref{eqn:StationaryLimitContra1} and \eqref{eqn:StationaryLimitContra2} contradict to each other. In only consequence, assumption \eqref{eqn:lem:StationaryLimit:WrongAssumption} must be wrong. q.e.d.
\end{itshape}

We admit that due to the small value of $\omega$ it can happen that $C_\phi$ becomes very large. This is why in Section~4.2 we include a practical enhancement that yields global convergence in a satisfactory amount of iterations regardless of the value that is chosen for $\omega$.

\subsection{Locally second order convergence}
From the Taylor series of the functions $\phi$ and $q$ at $\bx_k$ we find
\begin{align*}
	\nabla \phi(\bx_k + \alpha_k \cdot \bv_k) -\nabla q(\bx_k + \alpha_k \cdot \bv_k) =  \Big(\,\bH_k - \tbH_k\,\Big) \cdot \alpha_k \cdot \bv_k + \cO(\|\alpha_k \cdot \bv_k\|_2^2)\,.
\end{align*}
In the beginning of Section~2 we discussed sufficient conditions under which the choice $\tbH_k = \bH_k$ is suitable. The values for $\hbx_k = \bx_k + \bv_k$ then satisfy
\begin{align*}
	\nabla q(\hbx_k) = \bO
\end{align*}
because $\hbx_k$ is the unique minimizer of $q$. Thus, if $\alpha_k=1$ and $\tbH_k=\bH_k$ hold then
\begin{align*}
	\|\nabla \phi(\bx_k + \alpha_k \cdot \bv_k)\|_2 = \cO(\|\bv_k\|_2^2)\,.
\end{align*}

In the open neighborhood of a local minimizer $\bx^\star$ of $\phi$ it holds $\nabla \phi(\bx^\star)=\bO$ and $\nabla^2 \phi(\bx^\star)>0$. Thus $\phi$ can be approximated of second order by a parabola on the line $\bx_k(\alpha)=\bx_k + \alpha \cdot \bv_k$. In consequence of this, the line-search will return $\alpha_k$ sufficiently close to $1$ for all $k$ sufficiently large.

From the global convergence of $\lbrace\bx_k\rbrace$, discussed in Section~3.1, and the Cauchy criterion we find that there is an iteration $k$ from which holds $\|\bv_k\|_2 \ll 1$. We thus showed that for $k \in \N$ sufficiently large there will be second order convergence of
\begin{align}
	\|\nabla \phi(\bx_k)\|_2 \rightarrow 0\,.
\end{align}
We admit that the requirement on $\tbH$ is quite strong and may not hold for the problem instance at hand. This is why in Section~4.3 we provide an enhancement that guarantees second order local convergence under any circumstances.

\section{Practical enhancements}
Our algorithm uses four practical enhancements. These are:
\begin{itemize}
	\item a practical line-search;
	\item a watchdog technique \cite{watchdog} to accelerate global convergence;
	\item an additional Newton step per iteration to yield second-order convergence under no requirements;
	\item a primal-dual long-step interior-point method for solving the subproblems defined in \eqref{eqn:MinQ}.
\end{itemize}

\subsection{Line-search}
It is not worth the effort to use optimal values for $\alpha$. In practice it is important that the line-search makes the choice $\alpha=1$ under mild conditions, so that second-order local convergence can be easily achieved. On the other hand, it is important that also $\alpha>1$ can also be chosen because crude Hessian approximations for $\tbH$ can be overly convex, which results in very small length of the step direction $\|\bv\|_2$. In the following we propose a line-search that is cheap and accomplishes both goals.

The following line-search code shall replace line~\ref{algo:Solver:lineAlpha} in Algorithm~\ref{algo:Solver}. For ease of notation, we dropped the iteration index $k$ for $\bx_k,\bv_k,\alpha_k$.
\begin{algorithmic}[1]
\State $\alpha_{\max} := 1$\,,\quad $\alpha:=0$\,,\ $\cbx:=\bx$
\While{( \textbf{true} )}
	\State $\calpha:=\alpha_{\max}$
	\While{\textbf{not}$\Big($ \texttt{StepCriterion}($\cbx,\calpha,\bv$) \ \textbf{and} \  $\cbx+\calpha\cdot\bv \in \Omega$ $\Big)$}
		\State $\calpha:=\beta \cdot \calpha$
	\EndWhile
	\State $\cbx := \cbx + \calpha \cdot \bv$\,,\quad $\alpha := \alpha + \calpha$
	\If{( $\alpha<\alpha_{\max}$ )}
		\State \textbf{break}
	\EndIf
	\State $\alpha_{\max} := 2 \cdot \alpha_{\max}$
\EndWhile
\end{algorithmic}
The line-search comprises back-tracking with an iterative increase of the maximum trial step-size. But note, when $\alpha$ exceeds $\alpha_{\max}$ then we require that the step-criterion holds from the an updated value along the line. Thus, the step-criterion becomes more restrictive the more often we increase $\alpha_{\max}$. Typically, the step-criterion is Armijo's rule, i.e.
\begin{align*}
	\texttt{StepCriterion}(\bx,\alpha,\bv) := \phi(\bx + \alpha \cdot \bv) \leq \phi(\bx) + \gamma \cdot \alpha \cdot \nabla\phi(\bx)\t\cdot\bv\,.
\end{align*}
For the line-search parameters we choose $\gamma = 0.1$, $\beta = 0.8$, as suggested in \cite{Boyd}. However, we occasionally allow different step criteria, cf. Section~\ref{sec:Watchdog}\,.

\subsection{Watchdog}\label{sec:Watchdog}
The watchdog-technique is a particular line-search technique that is introduced in \cite{watchdog}. The motivation of watchdog is that for small values of $\omega$ and non-linear constraint functions $c$ the step-criterion due to Armijo will only allow very small steps because $\phi$ grows rapidly when $\|c(\bx)\|_2$ increases. Small step-sizes for $\alpha$ however mean that the algorithm would make little progress per iteration, resulting in large amounts of iterations and long computation times.

A way out of this dilemma is the use of a \textit{relaxed} step-criterion. The relaxed criterion admits larger values for $\alpha$ in the line-search and thus offers the convergence in a smaller amount of iterations, compared to the \textit{standard}, i.e. Armijo, step-criterion.

Implementing the algorithm with only a relaxed step-criterion is insufficient, as the relaxed condition is not restrictive enough to guarantee global convergence. The watchdog is an algorithmic safeguard that keeps track of the iterates. It tells our optimization method which step-criterion to use. The watchdog is aggressive, meaning that it would always allow our method to use the relaxed criterion, hoping it yields rapid convergence. But if the watchdog notices that the iterates won't make progress, it switches over to the standard criterion in order to force global convergence.

Strong theoretical results are available that prove that the watchdog-technique maintains the original global and local convergence properties of the algorithm. For all details on the implementation of the watchdog technique we refer the reader to \cite{watchdog}. We use a watchdog parameter $\ell=5$.

In our algorithm we use the following relaxed step-criterion. Call $\bx_\alpha := \bx + \alpha \cdot \bv$. For the relaxed step acceptance we require that at least either of these two conditions is satisfied.
\begin{align*}
	&\text{Condition 1} & &\phi(\bx_\alpha)<\phi(\bx)\\
	&\text{Condition 2} & &f(\bx_\alpha) < f(\bx)\quad \wedge \quad \|c(\bx_\alpha)\|_\infty \leq 10 \cdot \max\big\lbrace\,\|c(\bx)\|_\infty\,,\,0.01\,\big\rbrace
\end{align*}
The motivation for the above conditions is that we want a relaxed step-acceptance criterion while also avoiding totally unreasonable steps. The first condition says that there is progress after all. The second condition says that the objective improves while the constraint violation does not grow too much. Typically, when the constraint violation is moderately small then usually the subsequent descent steps will always be able to rapidly reach back to very small values for the norm of $c$, thus it is fine to use a maximum expression for further relax the criterion.

\subsection{An additional Newton step}
In Section~3.2 we proved second order local convergence of Algorithm~\ref{algo:Solver} to a stationary point under certain requirements. But in general, second-order convergence is impossible to achieve when only using step-directions that are obtained from solving \eqref{eqn:MinQ}. This is because $q$ is a convex approximation only, which can be insufficient. Letting $q$ be a nonconvex approximation could result in prohibitive cost for solving the subproblem \cite{Murty1987} and is thus not considered a practical option. But second order convergence can be achieved by using an additional search direction in Algorithm~\ref{algo:Solver} that arises from performing a simple Newton step.

This subsection is organized as follows. We first give a simple example problem that shows that steps obtained from the solution of a convex approximation $q$ do not permit second order convergence in general. We then the additional Newton step. Finally we discuss why indeed this step is sufficient to yield second-order convergence under any circumstances.

\paragraph{Example problem with no second-order convergence}
Consider the minimization of $f(x) = -0.5 \cdot x^2$ for $-1\leq x \leq 1$. The function $\phi$ in our algorithm becomes
\begin{align*}
	\phi(x) = (\rho-0.5) \cdot x^2 - \tau \cdot \Big( \log(1+x) + \log(1-x) \Big)
\end{align*}
$\omega$ does not appear since there are no equality constraints. Using the initial guess $x_0 = 0.5$, we hope to converge to a value close to $x=1$. For simplicity, we omit the convexization with $\rho$ and the left barrier term, yielding
\begin{align*}
	\phi(x) = -0.5 \cdot x^2 - \tau \cdot \log(1-x)\,.
\end{align*}
If for $q$ we use the positive semi-definite best-approximation $\tbH=0$ to $\bH=-1$, then $q$ has the following form at an iterate $x_k$:
\begin{align*}
	q_k(x) = -x_k \cdot x + \frac{\tau}{1-x}
\end{align*}
Using step-sizes of $\alpha=1$, as usually required for second-order convergence in higher dimensions, we get $x_{k+1} = \operatornamewithlimits{argmin}_{x \in \R} q_k(x)$. Since $q_k$ is convex, we can use the necessary optimality condition to obtain the explicit formula
\begin{align*}
	x_{k+1} = 1-\frac{\tau}{2 \cdot x_k}\,.
\end{align*}
This sequence converges to $x^\star = 0.5 + \sqrt{0.25 - 0.5 \cdot \tau}$ at a linear rate only, namely
\begin{align*}
	|x_{k}-x^\star| \in \Theta(\tau^k)\,.
\end{align*}

\paragraph{Second-order convergence of the Newton step}
We propose to add the following lines after line~\ref{algo:Solver:line:End} in Algorithm~\ref{algo:Solver}.
\begin{algorithmic}[1]
\State Attempt computing $\tbx := \bx_k - \nabla^2 \phi(\bx_k)\inv \cdot \nabla\phi(\bx_k)$\,.
\If{( $\tbx \in \R^n$ \ \textbf{and} \ $\phi(\tbx)<\phi(\bx_k)$)}
	\State $\bx_k := \tbx$
\EndIf
\end{algorithmic}
I.e., we attempt performing one Newton step for solving $\nabla\phi(\bx)=\bO$ from the initial guess $\bx_k$. This can fail since $\nabla^2 \phi(\bx_k)$ may be singular away from a local minimizer. But from above we know that in the local neighborhood of a minimizer it will be regular.

Since we search a local minimizer, we only accept the step if it yields reduction of the objective. Since \eqref{eqn:MinPhi} is unconstrained and strictly convex in the local minimizer (thanks to $\rho \cdot \lambda_\text{min}>0$), the above Newton step is second-order convergent whenever $\bx_k$ is sufficiently accurate. Since the sequence $\lbrace \bx_k \rbrace$ is globally convergent, eventually $\bx_k$ is sufficiently accurate.

\subsection{Primal-dual long-step interior-point method for solving the subproblems}
For sake of a self-contained presentation and for commenting on practical adaptations of this method, we discuss the algorithm introduced in \cite{StableIPM} that is used within our implementation of Algorithm~\ref{algo:Solver} for the solution of the subproblems \eqref{eqn:MinQ}.

We state the problem:
\begin{align*}
	\min_{\bx \in \Omega}\quad q(\bx) := &\frac{1}{2} \cdot \bx\t\cdot\bQ \cdot \bx + \bc\t\cdot\bx + \frac{1}{2 \cdot \omega} \cdot \|\bA \cdot \bx - \bb\|_2^2 \\
	&\quad - \tau \cdot \be\t \cdot \Big(\,\log(\be+\bx)+\log(\be-\bx)\,\Big)
\end{align*}
This is essentially a convex quadratic function augmented with barrier terms for box constraints. We define an auxiliary variable 
\begin{align}
	\blambda = \frac{-1}{\omega} \cdot (\bA \cdot \bx - \bb)\,.\label{eqn:def:lambda}
\end{align}

Our algorithm makes use of the following two functions, that are both parametric in $\nu>0$:
\begin{align*}
	\psi_\nu(\bx) := &\frac{1}{\nu} \cdot \Bigg(\, \frac{1}{2} \cdot \bx\t\cdot\bQ \cdot \bx + \bc\t\cdot\bx + \frac{1}{2 \cdot \omega} \cdot \|\bA \cdot \bx - \bb\|_2^2 \,\Bigg)\\
	& \quad -\Big(\, \log(\be+\bx)+\log(\be-\bx) \,\Big)\tageq\\[10pt]
	F_\nu(\bz) := & \begin{pmatrix}
		\bQ \cdot \bx + \bc - \bA\t \cdot \blambda - \bmu_L + \bmu_R\\
		\bA \cdot \bx - \bb + \omega \cdot \lambda\\
		\opdiag(\bmu_L) \cdot (\be+\bx)- \nu \cdot \be\\
		\opdiag(\bmu_R) \cdot (\be-\bx)- \nu \cdot \be
	\end{pmatrix}\,,\tageq
\end{align*}
where we use the short-hand $\bz = (\bx,\blambda,\bmu_L,\bmu_R) \in \R^{n+m+n+n}$.
\largeparbreak
\paragraph{Idea of the algorithm}
For all details on the algorithm we refer to \cite{StableIPM}. In the following we only give the ideas. Algorithm~\ref{algo:PrimalDual} states the method.

$\psi_\nu$ is self-concordant, strictly convex, and has a unique minimizer that converges to $\bx=\bO$ as $\nu\rightarrow +\infty$. One can prove that there is a value for $\nu$ that scales weakly with the logarithm of the norms of $\bQ,\bc,\bA,\bb$ and $\omega$ such that a Newton iteration for minimization of $\psi_\nu$ with initial guess $\bx=\bO$ converges rapidly to a sufficiently good minimizer of $\psi_\nu$ for the aforementioned suitable value of $\nu$. $\bx$ will then satisfy $2 \cdot \bx \in \Omega$. The suitable value for $\nu$ is found iteratively by evaluating an upper bound of the Newton-decrement at $\bx=\bO$, cf. \cite{Boyd,StableIPM} for details. 

In Algorithm~\ref{algo:PrimalDual}, the suitable value for $\nu$ is determined iteratively in line 4. Ihe sufficiently accurate minimizer of $\psi_\nu$ is computed with $10$ Newton iterations in line 7. Damping is not needed because due to the choice of $\nu$ it holds that $\bx$ is always sufficiently close to the exact minimizer so that the step-length $1$ in the Newton-iteration is always acceptable.

Once that the minimizer $\bx$ of $\psi_\nu$ is computed, we then augment the primal vector $\bx$ to a primal-dual vector $\bz$ by computing $\blambda$ as given above and
\begin{align}
	\bmu_L := \tau/(\be+\bx)\,,\quad \bmu_R := \tau/(\be-\bx)\,.\label{eqn:def:mu}
\end{align}
Since $\bx$ was an accurate root of $\nabla \psi_\nu$, it follows that $\bz$ is an accurate root of $F_\nu$ for the same value of $\nu$. Exact measures for the accuracy are given in \cite{StableIPM}. Given $\bz$ and $\nu$, we employ Mehrotra's predictor-corrector method to follow the path of roots of $F_\nu$ for iteratively reduced values of $\nu$ that converge to $\tau$. Eventually we arrive at a vector $\hbz$ satisfying
\begin{align*}
	F_\tau(\hbz)=\bO\,.
\end{align*}
The first component $\hbx$ of $\hbz$ solves our minimization problem \eqref{eqn:MinQ}. At least for the original short-step path-following version strong theory is available: All iterates $\bx$ are bounded away from $\partial\Omega$. All values in $\bmu_L,\bmu_R$ are bounded from below by strictly positive values. Last, the condition number of the Jacobian $DF_\nu$ of $F_\nu$ is bounded by a reasonable value at all iterates, even when the path-following iterates are perturbed by round-off errors \cite{StableIPM}.

For the long-step variant such guarantees do not exist. However, we find in practice that on averages it converges in $15$ iterations. The long-step path-following iteration uses the Mehrotra heuristic and is implemented in lines 11--23 in Algorithm~\ref{algo:PrimalDual}. Mehortra's method uses an affine step to estimate a value $\sigma$ for the geometric reduction of $\nu$, cf. line 17\,. A corrector step is then employed in order to restore centrality. Details on Mehrotra's method can be found in \cite{Mehrotra}.

Within the algorithm the set
\begin{align*}
	\cF := \Omega \times \R^m \times \R_+^n \times \R^n_+
\end{align*}
is used. The relation $\bz \in \cF$ means that its components $\bx$ are strictly interior and $\bmu_L,\bmu_R$ are strictly positive.

\begin{algorithm}
	\caption{Primal-dual method}
	\label{algo:PrimalDual}
	\begin{algorithmic}[1]
		\Procedure{CcpPbSolver}{\,$\bQ,\bc,\bA,\bb,\omega,\tau,\tol$\,}
		\State $\nu:=1$,\quad $\bx:=\bO$
		\While{$\|\nabla \psi_\nu(\bx)\|_2 \geq 0.25$}
			\State $\nu := 10 \cdot \nu$
		\EndWhile
		\For{$k=1,...,10$}
			\State $\bx := \bx - \nabla^2 \psi_\nu(\bx)^{-1} \cdot \nabla \psi_\nu(\bx)$
		\EndFor
		\State Compute $\blambda,\bmu_L,\bmu_R$ and state $\bz$, according to \eqref{eqn:def:lambda},\eqref{eqn:def:mu}.
		\While{$ \|F_\tau(\bz)\|_\infty >\tol $}
			\State $\nu := 0.5 \cdot \big(\,\bmu_L\t\cdot(\be+\bx) + \bmu_R\t\cdot(\be-\bx) \,\big)$
			\State \textit{// predictor (affine step to target $\nu=\tau$)}
			\State $\Delta\bz^{\text{aff}} := -DF_\nu(\bz)^{-1} \cdot F_\tau(\bz)$
			\State Choose $\alpha^{\text{aff}} \in (0,1]$ maximal subject to $\bz + \alpha^{\text{aff}} \cdot \Delta\bz^{\text{aff}} \in \overline{\cF}$
			\State $\bz^{\text{aff}} := \bz + \alpha^{\text{aff}} \cdot \Delta\bz^{\text{aff}}$
			\State $\nu^{\text{aff}} := 0.5 \cdot \big(\,(\bmu^{\text{aff}}_L)\t\cdot(\be+\bx^{\text{aff}}) + (\bmu^{\text{aff}}_R)\t\cdot(\be-\bx^{\text{aff}}) \,\big)$
			\State $\sigma := (\nu^{\text{aff}} / \nu)^3$
			\State $\hnu := \max\lbrace\,\tau\,,\,\sigma\cdot\nu\,\rbrace$
			\State \textit{// corrector (step to target $\nu=\hnu$)}
			\State $\Delta\bz^\text{cor} := -DF_\nu(\bz)^{-1} \cdot F_\hnu(\bz^{\text{aff}})$
			\State $\Delta\bz := \Delta\bz^{\text{aff}} + \Delta\bz^\text{cor}$
			\State Choose $\alpha \in (0,1]$ maximal subject to $\bz + \alpha \cdot \Delta\bz \in \overline{\cF}$
			\State $\bz := \bz + 0.99 \cdot \alpha \cdot \Delta\bz$
		\EndWhile
		\State \textit{// $\bz = (\bx,\blambda,\bmu_L,\bmu_R)$}
		\State \Return $\bx$
		\EndProcedure
	\end{algorithmic}
\end{algorithm}

\paragraph{Linear systems}
The linear systems to be solved in line 7 are well-posed since strongly dominated by the positive diagonal elements from the logarithmic terms in $\psi_\nu$. The linear systems in lines 13 and 20 can be made symmetric with the system matrix
\begin{align*}
	\bK = \begin{bmatrix}
	\bQ & \bA\t & \bI & \bI \\
	\bA & -\omega \cdot \bI & \bO & \bO\\
	\bI & \bO & - \opdiag\Big(\frac{\bmu_L}{\be+\bx}\Big) & \bO\\
	\bI & \bO & \bO & - \opdiag\Big(\frac{\bmu_L}{\be-\bx}\Big)
	\end{bmatrix}\,.
\end{align*}
It is symmetric indefinite, belonging to the category
\begin{align*}
	\bK= \begin{bmatrix}
		\bQ & \bG\t\\
		\bG & -\bD
	\end{bmatrix}\,,
\end{align*}
where $\bQ$ and $\bD$ are both positive definite. Thus,
\begin{align}
	\cond_2(\bK) \leq ( \|\bQ^{-1}\|_2 + \|\bD^{-1}\|_2 ) \cdot \|\bK\|_2 \,.
\end{align}
The inverse norm of $\bD$ in turn can be bounded from the lower and upper bounds that hold for all the iterates $\bx,\bmu_L,\bmu_R$. For details we refer to \cite{StableIPM}. The system can be easily reduced, using the Schur-complement
\begin{align}
	\Sigma := \bQ + \bG\t \cdot \bD^{-1} \cdot \bG\,,
\end{align}
where $\bD$ is a diagonal matrix.

\section{Numerical experiments}
For the numerical experiments we are particularly interested in large sparse NLP problems that arise from the discretization of one-dimensional path-constrained optimal-control problems. For the direct transcription of the control problem into a large space NLP we use the method introduced in \cite{StableTranscription}. This method yields an optimization problem where $\rho,\omega,\tau,\bS$ and functions for $f,c,\nabla f, \nabla c, \nabla^2_{\bx\bx} \cL$ as well as a sparsified positive semi-definite projection of $\nabla^2_{\bx\bx} \cL$ are provided. The objective that must be minimized is $\phi$ itself, so our method can be directly applied to solve these problems.

We compare our method against \Ipopt \cite{Ipopt}. Unfortunately, the interface of \Ipopt forbids to pass problems where $m$, the output-dimension of $c$, is larger than $n$. This is why we introduce auxiliary variables $\bs$ that we force by equality constraints to  satisfy $\omega \cdot \bs + c(\bx) = \bO$. Since $\|c(\bx)\|_2$ will be very small at the minimizer (in $\cO(\omega)$), it will be $\|\bs\|_2 \in \cO(1)$, i.e. the problem is reasonably scaled and no large numbers are introduced in the interface to \Ipopt. To use \Ipopt, we choose the objective $f(\bx) + \rho/2 \cdot \|\bx\|_\bS^2 + 0.5 \cdot \omega \cdot \|\bs\|^2$. Since \Ipopt uses primal barrier functions, the logarithmic terms with $\tau$ in $\phi$ will also appear in the effectively minimized objective of \Ipopt, yielding that both compared algorithms effectively solve the same optimization problem. For further details on how the problem is formulated to pass it to \Ipopt we refer to \cite[Section 5]{StableTranscription}.

The implementations are in \textsc R2016b in Windows 8.1 with Processor Intel(R) Core(TM) i7-4600U and 8GB RAM. We used the MEX-compiled \Ipopt version 3.11 from COIN-OR (www.coin-or.org).

\paragraph{Test problems}

Our test problems are listed in Table~\ref{table:Problems}. Problems 1 and 10 are given in \cite[eqns. 3 and 26]{Kameswaren}. Problems 2 and 3 can be found on the web-page of GPOPS-II (www.gpops2.com). The other problems are in order from \cite[pp. 79, 163, 85, 149, 113, 39]{BettsCollection}. Since we only access a computer with $\epsMach=10^{-16}$, we cannot choose $\rho,\omega,\tau$ very small. We solve all problems with $\rho=10^{-6}, \omega=10^{-6}, \tau=10^{-8}$. In all except two cases we use $\tol= 10^{-8}$. The exceptions are problem 5 and 7, since these are badly scaled.
\vspace{2mm}

The Aly-Chan problem is a problem with a singular arc where the sensitivity of the optimality value with respect to a variation in the control is below $10^{-10}$. Thus, it is difficult for the optimizer to find the unique smooth solution for this control, potentially resulting in many iterations.

The problems of brachistochrone, and those due to Bryson and Denham, Goddard, Hager and Rao form a biased set of well-scaled trial problems, involving convex quadratic programs and non-linear programs of small to large size.

Problems 5 and 7 involve biological models. As typical for biological problems, the states/controls differ widely in scales, leading to a bad scaling of the NLPs In our experiments we were unable to solve these NLPs to small tolerances.

Regarding the size, the problems 8 and 11 are most prohibitive. Problem 8 involves twelve species over a time-interval of 12 seconds. Problem 11 is originally stated for a time-interval of $10^4$ seconds. For the purpose of our experiments we reduced the time-interval to $100$ seconds, which still yields the same problem characteristic in terms of the shape of the solution.

For the initial guesses we used constant values for problems 1, 2, 3, 4, 8, 10 and 11. For problems 5, 6, 7 and 9 we used constant values for the controls and integrated the states numerically for the given initial conditions and constant control values. However, since most of the listed test problems involve end conditions, our initial guesses are usually infeasible.

The mesh-size has been chosen sufficiently small to yield curves for the discrete solutions that do qualitatively represent the shapes of the reference solutions. Yet, the mesh size is rather moderate, so that the NLPs are small enough to be solvable in reasonable amount of times on our computing system. Particularly the computation times of \Ipopt were a limiting factor in this regard.

\begin{table}
	\begin{tabular}{|l|l|l|l|l|l|}
		\hline
		\# & Name              & type                 & initial guess             & $h$      & FEM degree \\ \hline
		1  & Aly-Chan          & non-convex QCQP      & $\vec{0}$                 & $\pi/20$ & 10         \\ 
		2  & Brachistochrone   & NLP                  & $0.5 \cdot \vec{1}$       & $1/40$   & 4          \\ 
		3  & Bryson-Denham     & convex QP            & $\vec{0}$                 & $1/20$   & 8          \\ 
		4  & Chemical reactor  & NLP                  & $0.5 \cdot \vec{1}$       & $1/20$   & 4          \\ 
		5  & Chemotherapy      & badly scaled NLP     & $\int$, $\vec{u}=\vec{0}$ & $5$      & 4          \\ 
		6  & Caintainer crane  & NLP                  & $\int$, $\vec{u}=\vec{0}$ & $9/40$   & 4          \\ 
		7  & Drug treatment    & badly scaled NLP     & $\int$, $\vec{u}=\vec{0}$ & $5/2$    & 4          \\ 
		8  & Free flying robot & mildly nonlinear NLP & $\vec{1}$                 & $3/10$   & 8          \\ 
		9  & Goddard problem   & NLP     			  & $\int$, $T=T_m$           & $1/40$   & 4          \\ 
		10 & Hager problem     & convex QP            & $\vec{0}$                 & $1/4$    & 2          \\ 
		11 & Rao Problem       & NLP                  & $\vec{0}$                 & $1/10$   & 4          \\
		\hline
	\end{tabular}
	\caption[]{List of test problems: problem number, name, type of resulting finite-dimensional optimization problem, construction of initial guess, mesh-size, polynomial degree of finite element shape functions.}\label{table:Problems}
\end{table}

\paragraph{Experimental results} For the experiments we used a computation time limit of 10 hours and an iteration limit of 10000 iterations. \Ipopt caused two times a crash of \textsc{Matlab} for problems 2 and problem 11. In problem 11 it broke shortly before the time limit was reached. For the problems 2, 4, 5, 7 and 9 \Ipopt reached the iteration limit without a sufficiently accurate solution. For problem 1 \Ipopt terminated prematurely with a solution that it found of "acceptable accuracy", although $\tol=10^{-8}$ was specified. For problems 6 and 11 \Ipopt reported algorithmic errors in its own subroutines and terminated. In all except the aforementioned cases both algorithms worked as intended.

The results are given in Table~\ref{tab:ExperimentResults}. \Ipopt solved 2 out of 11 problems with success. PPD-IPM solved 11 out of 11 problems with success. All successfully solved problems of both solvers yielded feasible solutions that were good approximations to the reference optimal control solutions.

The table shows that in principle the iterations in \Ipopt are cheaper than in PPD-IPM. This is because \Ipopt solves only one linear system per iteration, while PPD-IPM utilizes Algorithm~\ref{algo:PrimalDual}, which solves $\approx 15$ linear systems per iteration. PPD-IPM compensates the larger cost per iteration by achieving a smaller amount of iterations in total. We further enhanced the solution of the linear systems in PPD-IPM by exploiting the particular sparsity pattern that results from the discretization. This explains why our method is only ten times slower per iteration, although it is \textsc{Matlab} (instead of MEX-compiled C++) and solves more systems per iteration.

We observe that for the convex quadratic programs in problems 3 and 10 PPD-IPM converged in one iteration. This is because the optimization of the subproblem $q$ is also an accurate interior-point solution to a convex quadratic program, confer to \cite{StableIPM}. Problem 9 converged in 4 iterations. During our investigation we found that our integrated initial guess with $T(t)=T_m$ is indeed the optimal solution, so the initial guess is -- apart from errors of numerical integration and discretization -- identical to the minimizer. PPD-IPM can strike a benefit of this accurate initial guess in that it converges in very few iterations, while \Ipopt cannot utilize the good initial guess.

It strikes the eye that in 7 out of 11 cases PPD-IPM terminates with a solution tolerance of less than $10^{-9}$ although the tolerance was only $\tol=10^{-8}$. This is because for each of the problems the method achieves local second order convergence. Thus, choosing smaller values for $\tol$ would only have negligible impact on the iteration count. Unfortunately, we cannot demonstrate this because $\epsMach=10^{-16}$ is not small enough.
\largeparbreak

In general, the results for \Ipopt are fairly bad, compared to PPD-IPM. We assume that the reason is as follows: We use exact Hessian matrices and we use highly accurate positive semi-definite approximations to the Hessian (instead of BFGS updates \cite[Chapter~6]{NumOpt}). PPD-IPM makes high use of the Hessian information. \Ipopt on the other hand does not. For example, when the Hessian is not positive definite, \Ipopt applies a crude shift with the identity, potentially resulting in bad progress because of the different scales in each solution component.

If instead of the exact Hessian we had used a crude Hessian approximation then certainly this would increase the iteration count of PPD-IPM. Then maybe the performance of \Ipopt would become better in comparison to PPD-IPM, since \Ipopt's iteration count would probably only change mildly while the count of PPD-IPM would grow severely. But we have accurate Hessians and want to strike a benefit of this, so there would be no point in doing experiments without using them.

\begin{table}
	\begin{tabular}{|l|l|l||l|l|l||l|l|l|}
		\hline
		\multicolumn{3}{c|}{Problem} & \multicolumn{3}{c|}{\Ipopt} & \multicolumn{3}{c|}{PPD-IPM}  \\
		\#      & n     & m     & time/iter    & \# iters    & NLP error              & time/iter & \# iters & $\|\nabla \phi\|_\infty$ \\ \hline\hline
		1       & 1505  & 1804  & 0.33         & 2188        & 9.8e-7$^{***,\dagger}$ & 1.3       & 146      & 6.2e-11                  \\ 
		2       & 2405  & 2885  & 0.30         & 10000       & 3.4e-2$^\dagger$       & 2.0       & 360      & 1.4e-9                   \\ 
		3       & 1443  & 1604  & 0.22         & 20          & 3.1e-9                 & 2.2       & 1        & 1.2e-11                  \\ 
		4       & 1924  & 2243  & 0.23         & 10000       & 3.3e-2$^\dagger$       & 1.5       & 138      & 3.5e-10                  \\ 
		5       & 6005  & 7204  & 0.65         & 10000       & 1.7e0$^\dagger$        & 2.6       & 275      & 7.4e-6                   \\ 
		6       & 3848  & 4492  & 0.63         & 1687        & 1.5e-1$^{*,\dagger}$   & 2.5       & 48       & 1.8e-9                   \\ 
		7       & 964   & 962   & 0.19         & 10000       & 1.0e+2$^{\dagger}$     & 1.4       & 203      & 2.6e-6                   \\ 
		8       & 11532 & 12812 & 20           & 543         & 9.5e-9                 & 12        & 71       & 5.1e-10                  \\ 
		9       & 2886  & 3525  & 0.30         & 10000       & 4.9e+1$^{\dagger}$     & 1.9       & 4        & 5.4e-10                  \\ 
		10      & 50    & 49    & 0.20         & 6           & 9.1e-10                & 1.4       & 1        & 1.7e-10                  \\ 
		11      & 24002 & 24002 & $\approx$8.5 & 4243$^{**}$ & ---                    & 5.9       & 13       & 1.2e-10                  \\
		\hline
	\end{tabular}
	\caption[]{Experimental results. Problem number, $n$ number of unknowns, each with box constraints, $m$ number of penalty-equality constraints; time per iteration, number of iterations for each solver. Regarding solution accuracy, \Ipopt monitors NLP error while we measure $\|\nabla\phi(\bx)\|_\infty$. Legend: $*$ Restoration failed; $**$ Restoration phase converged to feasible point that is unacceptable to the filter; $***$ \Ipopt says "acceptable solution"; $\dagger$ solution is not sufficiently accurate.}\label{tab:ExperimentResults}
\end{table}

\section{Final remarks}
We presented a novel optimization method that merges the ideas of primal interior-points, primal-dual interior-points, and successive quadratic programming. The method directly minimizes a penalty-barrier function, which is similar in approach to primal interior-point methods. In each iteration a step direction is computed by minimizing a convex subproblem that is the sum of a convex quadratic function, quadratic penalties, and logarithmic barriers. The approach of solving a convex (approximately quadratic) subproblem is related to successive quadratic programming. For the solution of the subproblem a primal-dual path-following method is used.

The method has some nice theoretical properties: It is very simple. There are no complicated issues related to infeasible subproblems because by construction all subproblems are feasible. The method has global convergence and second order local convergence (when an additional Newton step is used) under actually no requirements.

The method has three big practical advantages that make it yet un-competed by any other solver: First, it can solve overdetermined problems, i.e. where $m > n$, in a meaningful way in a comparably small amount of iterations. This is highly required, because the only convergent numerical scheme that is available for direct transcription of optimal control problems \cite{StableTranscription} results in such overdetermined problems. Also, the cost for solving the linear system is lower than the time for evaluating the problem, thus a small number of iterations is needed for good execution times. Second, our method does not require regularity of the Hessian of the Lagrange-function. Most other methods require it and thus enforce it by using a shift (that is significantly larger than $\rho$ and is not aligned with the objective), resulting in large amounts of iterations. Third, as the experiments showed our method can find highly accurate solutions to problems with indefinite or singular Hessian matrices and linearly dependent or overdetermined constraints. We believe this is only possible because highly accurate sparse approximations to the positive semi-definite projections of the exact Hessian matrices are employed for fast global convergence and exact Hessians are used for second-order local convergence. As far as we know, most other optimization methods do not even have an interface that would allow passing both a positive definite Hessian projection and the exact Hessian.

Further work shall be dedicated to implementing this method on a highly parallel system that computes with at least 32 significant digits. The algorithm itself is highly suitable for parallel computations because for our optimal control problems there are weakly scalable algorithms for the solution of the linear systems. From using 32 significant digits we expect that it removes all round-off related issues that could currently arise when choosing really small values for $\rho,\omega,\tau,\tol$.

\FloatBarrier

\bibliography{PPD_IPM_Bib}

\begin{thebibliography}{10}

\bibitem{BettsCollection}
J.T. Betts.
\newblock A collection of optimal control test problems.

\bibitem{Boyd}
Stephen Boyd and Lieven Vandenberghe.
\newblock {\em Convex Optimization}.
\newblock Cambridge University Press, New York, NY, USA, 2004.

\bibitem{watchdog}
R.~M. Chamberlain, M.~J.~D. Powell, C.~Lemarechal, and H.~C. Pedersen.
\newblock {\em The watchdog technique for forcing convergence in algorithms for
  constrained optimization}, pages 1--17.
\newblock Springer Berlin Heidelberg, Berlin, Heidelberg, 1982.

\bibitem{TrustRegionMethods}
Andrew~R. Conn, Nicholas I.~M. Gould, and Philippe~L. Toint.
\newblock {\em Trust-region Methods}.
\newblock Society for Industrial and Applied Mathematics, Philadelphia, PA,
  USA, 2000.

\bibitem{SUMT}
A.V. Fiacco and G.P. McCormick.
\newblock {\em Nonlinear Programming: Sequential unconstrained minimization
  techniques}.
\newblock John Wiley \& Sons.

\bibitem{Flet87}
Roger Fletcher.
\newblock {\em Practical Methods of Optimization}.
\newblock John Wiley \& Sons, New York, NY, USA, second edition, 1987.

\bibitem{Forsgren02interiormethods}
Anders Forsgren, Philip~E. Gill, and Margaret~H. Wright.
\newblock Interior methods for nonlinear optimization.
\newblock {\em SIAM Review}, 44:525--597, 2002.

\bibitem{Gill:2005:SSA:1055334.1055404}
Philip~E. Gill, Walter Murray, and Michael~A. Saunders.
\newblock {SNOPT}: An {SQP} algorithm for large-scale constrained optimization.
\newblock {\em SIAM Rev.}, 47(1):99--131, January 2005.

\bibitem{HOPDM}
Jacek Gondzio.
\newblock Hopdm (version 2.12) — a fast lp solver based on a primal-dual
  interior point method.
\newblock {\em European Journal of Operational Research}, 85(1):221 -- 225,
  1995.

\bibitem{IPM25ylater}
Jacek Gondzio.
\newblock Interior point methods 25 years later.
\newblock {\em European Journal of Operational Research}, 218(3):587--601,
  2012.

\bibitem{SQPinconsistent}
M.~T. Gouv$\hat{e}$a.
\newblock Dealing with inconsistent quadratic programs in a sqp based
  algorithm.
\newblock {\em Brazilian Journal of Chemical Engineering}, 14, 03 1997.

\bibitem{Han1977}
S.~P. Han.
\newblock A globally convergent method for nonlinear programming.
\newblock {\em Journal of Optimization Theory and Applications},
  22(3):297--309, Jul 1977.

\bibitem{Kameswaren}
Shivakumar Kameswaran and Lorenz~T. Biegler.
\newblock {Simultaneous dynamic optimization strategies: Recent advances and
  challenges}.
\newblock {\em Computers and Chemical Engineering}, 30(10):1560 -- 1575, 2006.
\newblock Papers form Chemical Process Control VII.

\bibitem{Kar84}
N.~Karmarkar.
\newblock A new polynomial-time algorithm for linear programming.
\newblock {\em Combinatorica}, 4(4):373--395, December 1984.

\bibitem{klee1970good}
V.~Klee and G.J. Minty.
\newblock {\em HOW GOOD IS THE SIMPLEX ALGORITHM}.
\newblock Defense Technical Information Center, 1970.

\bibitem{Mehrotra}
Sanjay Mehrotra.
\newblock {On the implementation of a primal-dual interior point method}.
\newblock {\em SIAM Journal on Optimization}, 2(4):575--601, 1992.

\bibitem{Murty1987}
Katta~G. Murty and Santosh~N. Kabadi.
\newblock Some np-complete problems in quadratic and nonlinear programming.
\newblock {\em Mathematical Programming}, 39(2):117--129, Jun 1987.

\bibitem{nesterov1994interior}
Y.~Nesterov and A.~Nemirovskii.
\newblock {\em Interior-point Polynomial Algorithms in Convex Programming}.
\newblock Studies in Applied Mathematics. Society for Industrial and Applied
  Mathematics (SIAM, 3600 Market Street, Floor 6, Philadelphia, PA 19104),
  1994.

\bibitem{StableTranscription}
M.~P. {Neuenhofen}.
\newblock {High-order convergent Finite-Elements Direct Transcription Method
  for Constrained Optimal Control Problems}.
\newblock {\em ArXiv e-print 1712.07761}, December 2017.

\bibitem{StableIPM}
M.~P. {Neuenhofen}.
\newblock {Stable interior-point method for convex quadratic programming with
  strict error bounds}.
\newblock {\em ArXiv e-print 1711.01418}, November 2017.

\bibitem{NumOpt}
Jorge Nocedal and Stephen~J. Wright.
\newblock {\em Numerical Optimization}.
\newblock Springer, New York, NY, USA, second edition, 2006.

\bibitem{Schittkowski1982}
Klaus Schittkowski.
\newblock The nonlinear programming method of wilson, han, and powell with an
  augmented lagrangian type line search function.
\newblock {\em Numerische Mathematik}, 38(1):83--114, Feb 1982.

\bibitem{LOQO}
Robert~J. Vanderbei.
\newblock Loqo:an interior point code for quadratic programming.
\newblock {\em Optimization Methods and Software}, 11(1-4):451--484, 1999.

\bibitem{IPOPT}
Andreas W{\"a}chter and Lorenz~T. Biegler.
\newblock On the implementation of an interior-point filter line-search
  algorithm for large-scale nonlinear programming.
\newblock {\em Mathematical Programming}, 106(1):25--57, Mar 2006.

\bibitem{Wright}
S.J. Wright.
\newblock {\em Primal-Dual Interior-Point Methods}.
\newblock Other Titles in Applied Mathematics. Society for Industrial and
  Applied Mathematics, 1997.

\bibitem{YildirimW02}
E.~Alper Yildirim and Stephen~J. Wright.
\newblock Warm-start strategies in interior-point methods for linear
  programming.
\newblock {\em {SIAM} Journal on Optimization}, 12(3):782--810, 2002.

\end{thebibliography}
\bibliographystyle{plain}

\end{document}